\documentclass{amsart}
\usepackage{amsmath}
\usepackage{amssymb}
\usepackage{latexsym}
\usepackage{amsfonts}

\newtheorem{Pa}{Paper}[section]
\newtheorem{Tm}[Pa]{{\bf Theorem}}
\newtheorem{La}[Pa]{{\bf Lemma}}

\newtheorem{Pb}[Pa]{{\bf Problem}}

\newcommand{\CC}{{\mathchoice
{\setbox0=\hbox{$\displaystyle\rm C$}\hbox{\hbox
to0pt{\kern0.4\wd0\vrule height0.9\ht0\hss}\box0}}
{\setbox0=\hbox{$\textstyle\rm C$}\hbox{\hbox
to0pt{\kern0.4\wd0\vrule height0.9\ht0\hss}\box0}}
{\setbox0=\hbox{$\scriptstyle\rm C$}\hbox{\hbox
to0pt{\kern0.4\wd0\vrule height0.9\ht0\hss}\box0}}
{\setbox0=\hbox{$\scriptscriptstyle\rm C$}\hbox{\hbox
to0pt{\kern0.4\wd0\vrule height0.9\ht0\hss}\box0}}}}

\begin{document}
\title{Rigidity, boundary interpolation and reproducing kernels}

\author[D. Alpay]{Daniel Alpay}
\address{(DA) Department of Mathematics\\
Ben-Gurion University of the Negev\\
POB 653\\ 84105 Beer-Sheva, Israel}
\email{dany@math.bgu.ac.il}
\author[S. Reich]{Simeon Reich}
\address{(SR) Department of Mathematics\\
      The Technion - Israel Institute of Technology\\
      32000 Haifa,
      Israel}
\email{sreich@techunix.technion.ac.il}
\author[D. Shoikhet]{David Shoikhet}
\address{(DS) Department of Mathematics\\
Ort Braude College\\
POB 78, 21982 Karmiel, Israel}
\email{davs27@netvision.net.il}

\begin{abstract}
We use reproducing kernel methods to study various rigidity
problems. The methods and setting allow us to also consider the
non-positive case.
\end{abstract}
  \maketitle
\section{Introduction}
\setcounter{equation}{0}
We use reproducing kernel methods to study various rigidity
problems. First, recall that a function $s$ is analytic and
contractive in the open unit disk if and only if the kernel
\[
k_s(z,w)=\dfrac{1-s(z)\overline{s(w)}}{1-z\overline{w}}
\]
is positive in the open unit disk ${\mathbb D}$. Such functions
have various names in the literature. We will use the term
{\it Schur functions} in the present paper, and will denote by
${\mathcal S}_0$ the family of Schur functions. We note that
constants in the closed unit disk, and, in particular,
unimodular constants are Schur functions. See \cite{Alpay01}
for a review of Schur functions. Krein and Langer
introduced in \cite{kl1} the class of {\it generalized Schur functions}
as the set of functions $s$ meromorphic in the open unit disk and
such that the kernel
\[
k_s(z,w)=\dfrac{1-s(z)\overline{s(w)}}{1-z\overline{w}}
\]
has a finite number (say, $\kappa$) of negative squares for $z$
and $w$ in the domain of analyticity $\Omega(s)$ of $s$. This
means that for every choice of an integer $n$ and of points
$w_1,\ldots, w_n$ in $\Omega(s)$, the $n\times n$ hermitian matrix
the $(\ell,j)$-th entry of which is equal to $k_s(w_\ell, w_j)$,
has at most
$\kappa$ strictly negative eigenvalues, and exactly $\kappa$
strictly negative eigenvalues for some choice of $n,w_1,\ldots ,
w_n$.\\

We denote by ${\mathcal S}_{\kappa}$ the family of {\it generalized
Schur functions} with $\kappa$ negative squares, and set
\[{\rm sq}_-(s)=\kappa.\]
By a result of Krein and Langer \cite{kl1}, a function
$s\in{\mathcal S}_\kappa$ if and only if it can be written as
\[
s(z)=\dfrac{s_0(z)}{b(z)},
\]
where $s_0\in{\mathcal S}_0$ and $b$ is a Blaschke product of
order $\kappa$, and where, moreover, $s_0$ and $b$ have no common
zeros. For example, the function $s(z)=1/z$ belongs to ${\mathcal
S}_1$.
 See, for instance, \cite{adrs} for
more information on generalized Schur functions. It is well to
recall at this point that the positivity of the kernel $k_s(z,w)$
implies the analyticity of $s$. On the other hand, a function $s$
can be such that the kernel $k_s(z,w)$ has a finite number of
negative squares in the open unit disk without being meromorphic
there, as is illustrated by the example
\[s(z)=\begin{cases}0,\quad z\not= 0,\\
                    1,\quad z=0.

                    \end{cases}
                    \]
The corresponding kernel $k_s(z,w)$ has one negative square, but
$s$ is not meromorphic in ${\mathbb D}$. See \cite[p. 82]{adrs}.
Thus, meromorphicity is an essential part of the definition of a
generalized Schur function.\\

The methods and setting of the present paper allow us to also
study rigidity problems in the non-positive case, where Schur
functions are replaced by generalized Schur functions.\\

We now turn to the outline of the paper. The paper consists of six
sections, including this introduction. In Section 2, we review
the main results on boundary interpolation which we will need.
Section 3 contains Theorem \ref{Malakoff}, which is the main
result of this paper. It is, to the best of our knowledge, the
first rigidity theorem where meromorphic functions are considered.
We illustrate Theorem \ref{Malakoff} in two cases, in Sections 4
and 5: in Section 4 we recover, in particular, a well-known theorem of
Burns and Krantz, while in Section 5 we consider an example in the
non-positive case. Finally, in Section 6, we consider a case
where, in general, there is no rigidity result.

\section{Boundary interpolation problem: a quick review}
\setcounter{equation}{0}
We review in this section some of the results proved in
\cite{ADLW} and which we will use in the sequel. We denote by
${\mathbb T}$ the unit circle and by $z \hat\rightarrow z_1$ the
nontangential (or angular) convergence of $z$ to $z_1\in{\mathbb T}$. This
means (see, for instance, \cite[p. 47]{sarason94}) that $z$ stays
in a region of the form
\begin{equation}
\label{julia}
\{ z \in \mathbb{D}: |z-z_1| <  K(|1-|z|)\}
\end{equation}
for some $K>1$.
\begin{Pb}
\label{pb1}
Let $z_1\in{\mathbb T}$, an integer $k\geq 1$, and
complex numbers $\tau_0$, $\tau_k$, $\tau_{k+1}, \dots ,
\tau_{2k-1}$ with $|\tau_0|=1$,  $\tau_k\ne 0$, be given. Find all
generalized Schur functions $s$ such that
\begin{equation}\label{exp}
s(z)=\tau_0+\sum_{i=k}^{2k-1}\tau_i(z-z_1)^i + {\rm
O}((z-z_1)^{2k}), \quad z \hat\rightarrow z_1.
\end{equation}
\end{Pb}

Before presenting the solution to this problem we introduce some
notations. First, we define the matrices
\begin{equation}\label{Bk}
 T=\begin{pmatrix}
  \tau_k& 0 & \cdots&0 &0\\
\tau_{k+1} & \tau_k & \dots  &0&0  \\
\vdots  &    \vdots  & \ddots&\vdots & \vdots \\
\tau_{2k-2}&\tau_{2k-3}&\cdots&\tau_k&0\\
\tau_{2k-1} &\tau_{2k-2}& \cdots &\tau_{k+1}& \tau_k
\end{pmatrix}
\end{equation}
and
\begin{equation}\label{bbb}
B=\begin{pmatrix}
0&0&\cdots&              0  &  (-1)^{k-1}\binom{k-1}{0} z_1^{2k-1} \\[2mm]
0     & 0  &\cdots        &(-1)^{k-2}\binom{k-2}{0} z_1^{2k-3}
&  (-1)^{k-1}\binom{k-1}{1} z_1^{2k-2}     \\[2mm]
\vdots  & \vdots          & \vdots        &    \vdots        &   \vdots      \\[2mm]
0      &-\binom10 z_1^3 &\cdots &            (-1)^{k-2}
\binom{k-2}{k-3} z_1^{k} &(-1)^{k-1}\binom{k-1}{k-2} z_1^{k+1}\\[2mm]
 z_1 &-\binom11 z_1^2& \cdots&(-1)^{k-2}\binom{k-2}{k-2} z_1^{k-1}
 &(-1)^{k-1}\binom{k-1}{k-1}
 z_1^k
\end{pmatrix}.
\end{equation}
We note that $B$ is a right lower triangular matrix, which is
invertible because $z_1\not=0$, and since
$\tau_k,\,\tau_0\not=0$, the matrix
\begin{equation}
\label{P} \mathbb P:=\overline{\tau_0}TB
\end{equation}
is also invertible. We also define the vector function
$$
 R(z)=\begin{pmatrix} \dfrac{1}{1-z\overline{z_1}} & \dfrac{z}{(1-z\overline{z_1})^2}
 & \ldots
& \dfrac{z^{k-1}}{(1-z\overline{z_1})^k}\end{pmatrix},
$$
fix some  $z_0\in \mathbb T$, $z_0\not =z_1$, and define the
polynomial
\begin{equation}
\label{polynomial} p(z)=(1-z\overline{z_1})^k \,R(z)\mathbb P^{-1}
R(z_0)^*.
\end{equation}
Its degree is at most $k-1$ and $p(z_1)\not =0$. See \cite[p.
13]{ADLW} for this last fact.

\begin{Tm}
\label{bbip}
Let
 $z_1 \in \mathbb T$ and $\tau_0, \tau_k, \ldots , \tau_{2k-1}$ be
as in Problem  $\ref{pb1}$, assume that the matrix $\mathbb P$ in
\eqref{P} is hermitian, and let $\Theta$ be the $J$--unitary
rational matrix function
\begin{equation}
\label{theta} \Theta(z)= \begin{pmatrix}
a(z)&b(z)\\c(z)&d(z)\end{pmatrix}= I_2-
\dfrac{(1-z\overline{z_0})p(z)}{(1-z\overline{z_1})^k}
{\mathbf u}
{\mathbf u}^*J
,\quad J=\begin{pmatrix}1&0\\0&-1\end{pmatrix},
\quad \mathbf u=
\begin{pmatrix} 1 \\ \overline{\tau_0} \end{pmatrix},
\end{equation}
with $p$ defined by \eqref{polynomial} and a fixed $z_0 \in
\mathbb T$, $z_0\neq z_1$. Then
 the fractional-linear transformation
\begin{equation}
\label{lft} s(z)=T_{\Theta(z)}(s_1(z))
=\dfrac{a(z)s_1(z)+b(z)}{c(z)s_1(z)+d(z)}
\end{equation}
establishes a bijective correspondence between all solutions $s$
of  Problem $\ref{pb1}$ and all
$s_1\in\cup_{\kappa=0}^\infty{\mathcal S}_\kappa$ with the
property
\begin{equation}
\label{liminf}
\liminf _{z \hat{\rightarrow}
z_1}|s_1(z)-\tau_0|>0.
\end{equation}
Moreover, if $s$ and $s_1$ are related by \eqref{lft}, then
\begin{equation}  \label{eqq}
{\rm sq}_-(s)={\rm sq}_-(s_1)+{\rm ev}_-({\mathbb P}),
\end{equation}
where ${\rm ev}_-({\mathbb P})$ denotes the number of strictly
negative eigenvalues of ${\mathbb P}$.
\end{Tm}

This is \cite[Theorem 3.2, p. 43]{ADLW}. The proof there is based
on the theory of reproducing kernel Pontryagin spaces of the kind
introduced by L. de Branges and J. Rovnyak. Having in mind the
proofs of the forthcoming theorems, it is well to recall that one
can write
\[
\Theta(z)=\begin{pmatrix}1-\theta(z)&\tau_0\theta(z)\\
-\overline{\tau_0}\theta(z)&1+\theta(z)\end{pmatrix},
\]
where
\[
\theta(z)=\frac{(1-z\overline{z_0})p(z)}{(1-z\overline{z_1})^k}.\]
See \cite[(3.18), p. 16]{ADLW}.\\

We have, in particular (see the formula on the sixth line from the
top of page 16 of \cite{ADLW}),
\[
c(z)s_1(z)+d(z)=-\overline{\tau_0}(s_1(z)-\tau_0)\theta(z)+1,
\]
and hence
\begin{equation}
\label{opera_bastille}
%
c(z)s_1(z)+d(z)
=\frac{(1-z\overline{z_1})^k-\overline{\tau_0}(1-z\overline{z_0})p(z)(s_1(z)-\tau_0)}
{(1-z\overline{z_1})^k}.
\end{equation}
This expression also makes clear why property \eqref{liminf} is
important. When combined with the fact that $p(z_1)\not=0$, it
yields the conclusion that the numerator in
\eqref{opera_bastille} stays bounded away from $0$ as $z$ tends
nontangentially to $z_1$.\\

Finally, we recall that the description of all the solutions of the
interpolation problem is independent of the choice of
$z_0\in{\mathbb T}\setminus\{z_1\}$. This follows from Theorem
\ref{bbip} since the problem itself does not depend on $z_0$.
For the convenience of the reader we present a direct argument. Let
$\Theta$ be the matrix function corresponding to $z_0$ and
$\widehat{\Theta}$ the matrix function corresponding to a
different point $\widehat{z_0}$ on the unit circle. Then, see
\cite[(3.22), p. 17]{ADLW}, the corresponding functions $\theta$
and $\widehat{\theta}$ are related by
\[
\theta(z)-\widehat{\theta}(z)=-
\widehat{\theta}(z_0).
\]
It follows that
\[
I_2-\theta(z){\mathbf u} {\mathbf u}^*J=\left(
I_2-\widehat{\theta}(z){\mathbf u} {\mathbf u}^*J\right) \left(
I_2+\widehat{\theta}(z_0){\mathbf u} {\mathbf u}^*J\right),\]
that is,
\begin{equation}
\label{la_defense}
\Theta(z)=\widehat{\Theta}(z)U,
\end{equation}
where $U$ is the $J$--unitary matrix
\[
U=I_2+\widehat{\theta}(z_0){\mathbf u} {\mathbf u}^*J=
\begin{pmatrix}1+\widehat{\theta}(z_0)&-\tau_0\widehat{\theta}(z_0)\\
\overline{\tau_0}\widehat{\theta}(z_0)&1-\widehat{\theta}(z_0)\end{pmatrix}.
\]
Let us use the notation
\[T_{\Theta(z)}(s_1(z))=\frac{a(z)s_1(z)+b(z)}{c(z)s_1(z)+d(z)}.\]
From \eqref{la_defense} it follows that
\[
T_{\Theta(z)}(s_1(z))=T_{\widehat{\Theta}(z)}\left(T_U(s_1(z))\right).
\]
Hence, to show that the description of the set of solutions does not
depend on the normalization point, we have to show that condition
\eqref{liminf} is invariant under the linear-fractional
transformation defined by $U$, that is, that the function
\[\widehat{s_1}(z)=
\frac{(1+\widehat{\theta}(z_0))s_1(z)-\tau_0\widehat{\theta}(z_0)}
{\overline{\tau_0}\widehat{\theta}(z_0)s_1(z)+1-\widehat{\theta}(z_0)}
\]
satisfies \eqref{liminf} if and only if $s_1$ does. But this is
clear from the equality
\[
\widehat{s_1}(z)-\tau_0=\frac{s_1(z)-\tau_0}{1+\widehat{\theta}(z_0)\overline{\tau_0}
(s_1(z)-\tau_0)}.
\]
\section{A general rigidity theorem}
\setcounter{equation}{0}
As we have already mentioned,  the theorem to be proved in this section
is, to the best of our knowledge, the first rigidity theorem where
meromorphic functions are considered. Let us first explain in
words this theorem. Consider the boundary interpolation problem
\ref{pb1}, and fix $z_1=1$ to lighten the notation. One shows
that any two solutions of this problem, $s$ and $\sigma$, satisfy
\[
s(z)-\sigma(z)=((z-1)^{2k})g(z),
\]
where the function $g$ stays bounded as $z$ tends nontangentially
to $1$. Using Julia's lemma, we then show that for a special
choice of $\sigma$, and if $s$ and $\sigma$ satisfy the stronger
requirement
\[
s(z)-\sigma(z)=O((z-1)^{2k+2}),
\]
then they must coincide: $s=\sigma$.
\begin{Tm}
\label{Malakoff}
Consider the boundary interpolation problem
\ref{pb1} with $z_1=1$, and assume that the associated matrix
${\mathbb P}$ is hermitian. Let $x\in{\mathbb
T}\setminus\left\{\tau_0\right\}$ and let
$b(z)=T_{\Theta(z)}(x)$. Let $s\in{\mathcal S}_{{\rm
ev}_-({\mathbb P})}$ be a solution of Problem \ref{pb1} such that
\begin{equation}
\label{Luxembourg} s(z)-b(z)=O((z-1)^{2k+2}).
\end{equation}
Then
\[
s(z)\equiv b(z).
\]
\end{Tm}

Note that $x\in{\mathbb T}\setminus\left\{\tau_0\right\}$ obviously
has property \eqref{liminf}.\\

 Theorem
\ref{Malakoff} includes a number of results, and, in particular,
the well-known Burns-Krantz rigidity theorem. We note that $b$ is
a rational function which takes unitary values on the unit
circle, that is, $b$ is a quotient of two finite Blaschke products.\\

Before proving Theorem \ref{Malakoff}, we present a preliminary lemma.

\begin{La}
\label{Richard-Lenoir}
Let $x\in{\mathbb T}$ and let $\sigma$ be a Schur function such that
\begin{equation}
\label{220707}
\sigma(z)-x=O((z-1)^2).
\end{equation}
Then $\sigma\equiv x$.
\end{La}

{\bf Proof:} Indeed, by \eqref{220707}, we have that
\[
\liminf_{z \hat\rightarrow 1}
\Big|\dfrac{\sigma(z)-x}{1-z}\Big|=0,\]
and hence
\[
c\,:\,\,=\liminf_{z \hat\rightarrow 1} \dfrac{1-|\sigma(z)|}{1-|z|}=0,\]
when $z$ is restricted to a region of the form \eqref{julia}.
Assume now by contradiction that $\sigma(z)\not \equiv x$. It then follows
from Julia's lemma (see, for instance, \cite[p. 51]{sarason94}) that
\[
\dfrac{|\sigma(z)-x|^2}{1-|\sigma(z)|^2}\le
c\dfrac{|z-1|^2}{1-|z|^2}=0
\]
for any $z\in{\mathbb D}$. This contradicts the hypothesis
$\sigma(z)\not \equiv x$.\mbox{}\qed\mbox{}\\

{\bf Proof of Theorem \ref{Malakoff}:} Using the formula
\[
\det (I_p-AB)=\det(I_q-BA),
\]
where $A\in{\mathbb C}^{p\times q}$ and $B\in{\mathbb C}^{q\times
p}$, and since
\[
{\mathbf u}^*J{\mathbf u}=0,
\]
we obtain from \eqref{theta} that
\[
\det
\Theta(z)=1-\dfrac{(1-z\overline{z_0})p(z)}{(1-z\overline{z_1})^k}{\mathbf
u}^*J{\mathbf u}=1.
\]
Let $s$ be a solution to Problem \ref{pb1}, which belongs to
${\mathcal S}_{\nu_-({\mathbb P})}$. There is a Schur function
$s_1$ subject to \eqref{liminf} such that $s=T_\Theta(s_1)$. Thus,
with $b(z)=T_{\Theta(z)}(x)$, we have
\[
\begin{split}
s(z)-b(z)&=\dfrac{(\det\Theta(z))(s_1(z)-x)}{(c(z)s_1(z)+d(z))(c(z)x+d(z))}\\
&=\dfrac{s_1(z)-x}{(c(z)s_1(z)+d(z))(c(z)x+d(z))}.
\end{split}
\]
Using \eqref{opera_bastille} for $s_1$ and for $x$, we have
\[
s(z)-b(z)=\dfrac{(1-z)^{2k}(s_1(z)-x)}{
((1-z)^k-\overline{\tau_0}(1-z\overline{z_0})p(z)(s_1(z)-\tau_0))
((1-z)^k-\overline{\tau_0}(1-z\overline{z_0})p(z)(x-\tau_0))}.
\]
Now condition \eqref{liminf} and the fact that $p(z_1)\not =0$
come into play. If \eqref{Luxembourg} holds, then it follows that
\begin{equation*}
s_1(z)-x=O((z-1)^2),
\end{equation*}
and hence $s_1(z)\equiv x$ by Lemma \ref{Richard-Lenoir}.\mbox{}\qed\mbox{}\\

\section{Recovering a theorem of Burns and Krantz}
\setcounter{equation}{0}
\label{lala}
Theorem \ref{Malakoff} specialized to  $k=1, \;
\tau_0=\tau_1=1$, and $x=-1$ reduces to Theorem
\ref{Chatelet} presented in this section. When the
angular convergence is replaced by the unrestricted one,
Theorem \ref{Chatelet} is due to D. M. Burns and S. G.
Krantz \cite{MR1242454}. For a survey of related results, see
\cite{ELRS}.

\begin{Tm}
Assume that a Schur function $s$ satisfies
\begin{equation}
s(z)=z+O((1-z)^4),\quad z\hat{\rightarrow} 1,
\label{Crimee, ligne 7}
\end{equation}
where $\hat{\rightarrow}$ denotes nontangential convergence. Then
$s(z)\equiv
z$.
\label{Chatelet}
\end{Tm}
We can recover this result as a special case of Theorem
\ref{Malakoff}. We present the argument for completeness. To
begin with, $s$ is, in particular, a solution of the boundary
interpolation problem
\begin{equation}
\label{voltaire} s(z)=z+O((1-z)^2),\quad z\hat{\rightarrow} 1.
\end{equation}
To solve this last problem we use Theorem \ref{bbip}.
In the notation of the theorem, we have
\[
z_1=1,\quad k=1\quad{\rm and}\quad\tau_0=\tau_1=1.
\]
The matrix ${\mathbb P}$ given by \eqref{P} reduces
to a number: ${\mathbb P}=\overline{\tau_0}\tau_1z_1=1$.\\

To compute the coefficient matrix--function (which will allow us
to describe the set of all solutions to \eqref{voltaire}), we need
to fix a point $z_0\in{\mathbb T}\setminus\{1\}$. We will choose
$z_0=-1$. With this choice of $z_0$,
the polynomial $p$ given in \cite[(3.2), p. 11]{ADLW}, the
definition of which is recalled in \eqref{polynomial}, is then:
\[
p(z)=(1-z)\frac{1}{1-z}1^{-1}\frac{1}{2}=\frac{1}{2}.
\]
The matrix function $\Theta$ in Theorem \ref{bbip} is then equal to
\begin{equation}
\begin{split}
\Theta(z)&=\begin{pmatrix}a(z)&b(z)\\
c(z)&d(z)\end{pmatrix}
=\begin{pmatrix}1&0\\0&1\end{pmatrix}-\frac{1+z}{2(1-z)}\begin{pmatrix}1\\1\end{pmatrix}
\begin{pmatrix}1&1\end{pmatrix}\begin{pmatrix}1&0\\0&-1\end{pmatrix}\\
&=\frac{1}{2(1-z)}\begin{pmatrix}1-3z&1+z \\-1-z&3-z
\end{pmatrix}.
\end{split}
\end{equation}
Thus, Theorem \ref{bbip} becomes here:

\begin{La}
The fractional-linear transformation
\begin{equation}
\label{lumbago1}
s(z)=\dfrac{(1-3z)s_1(z)+1+z}{3-z-(1+z)s_1(z)}
\end{equation}
gives a one-to-one correspondence between Schur functions
satisfying \eqref{voltaire},
\[
s(z)=z+O((1-z)^2),\quad z\hat{\rightarrow} 1,
\]
and Schur functions $s_1$ which satisfy condition \eqref{liminf}.
\end{La}
The choice $x=-1$ corresponds to $s(z)=z$.\\

We now compute $s(z)-z$:
\begin{equation}
\begin{split}
s(z)-z&=\dfrac{(1-3z)s_1(z)+1+z}{3-z-(1+z)s_1(z)}-z\\
&=\dfrac{(1-3z)s_1(z)+1+z-3z+z^2+(z^2+z)s_1(z)}{3-z-(1+z)s_1(z)}\\
&=(1-z)^2\dfrac{s_1(z)+1}{3-z-(1+z)s_1(z)}.
\end{split}
\label{Corentin Cariou, ligne 7}
\end{equation}
We intend to show that condition \eqref{Crimee, ligne 7} forces
$s_1(z)\equiv -1$. It will then follow from
\eqref{Corentin Cariou, ligne 7}
that $s(z)\equiv z$.\\

Combining \eqref{Corentin Cariou, ligne 7} and
\eqref{Crimee, ligne 7},
we obtain
\[
\dfrac{s_1(z)+1}{3-z-(1+z)s_1(z)}=O((1-z)^2),\quad
z\hat{\rightarrow} 1.
\]
But $|3-z-(1+z)s_1(z)|\le 6$ in the open unit disk, and therefore
\[
|s_1(z)+1|=|3-z-(1+z)s_1(z)|\cdot\left|\dfrac{s_1(z)+1}{3-z-(1+z)s_1(z)}\right|\le
K|1-z|^2,\quad z\hat{\rightarrow} 1,\]
for some positive constant $K$. Thus
\[
s_1(z)+1=O((1-z)^2),\quad z\hat{\rightarrow} 1.\]
By Lemma \ref{Richard-Lenoir}, $s_1(z)+1\equiv 0$, and this ends the proof.

\section{An example in the indefinite case}
\setcounter{equation}{0}
\begin{Tm}
Let $s$ be a generalized Schur function with one negative square
and assume that
\[
s(z)-\dfrac{1}{z}=O((1-z)^4).
\]
Then
\[
s(z)\equiv\dfrac{1}{z}.
\]
\end{Tm}
{\bf Proof:} The function $s$ is, in particular, a solution to
Problem \ref{pb1} with
\begin{equation}
\label{lumbago11}
 z_1=1,\quad k=1\quad{\rm
and}\quad\tau_0=1,\quad \tau_1=-1.
\end{equation}
The matrix ${\mathbb P}$  reduces
to a strictly negative number:
${\mathbb P}=\overline{\tau_0}\tau_1z_1=-1$. Thus
there are solutions to the interpolation problem in the class of
generalized Schur functions ${\mathcal S}_1$.
The polynomial $p$ is now equal to
\[
p(z)=(1-z)\frac{1}{1-z}\frac{1}{(-1)}\frac{1}{2}=-\frac{1}{2},\]
and we have, with $z_0=-1$,
\[
\Theta(z)=\begin{pmatrix}1&0\\0&1\end{pmatrix}+\frac{1}{2}
\frac{1+z}{1-z}\begin{pmatrix}1\\
1\end{pmatrix}\begin{pmatrix}1&-1\end{pmatrix}.\]
Thus
\[
2(1-z)\Theta(z)=\begin{pmatrix}2(1-z)+(1+z)&-1-z\\
1+z&2(1-z)-1-z\end{pmatrix} =
\begin{pmatrix}
3-z&-1-z\\1+z&1-3z\end{pmatrix}.
\]
By Theorem \ref{bbip}, a generalized Schur function $s$
with one negative square
satisfies the interpolation problem with data \eqref{lumbago11} if
and only if it is of the form
\begin{equation}
\label{lumbago33} s(z)=\dfrac{(3-z)s_1(z)-1-z}{1-3z+(1+z)s_1(z)},
\end{equation}
where the parameter $s_1$ is any
Schur function satisfying \eqref{liminf}.\\

We now compute the difference $s(z)-\dfrac{1}{z}$:

\[
\begin{split}
s(z)-\dfrac{1}{z}&=\dfrac{(3-z)s_1(z)-1-z}{1-3z+(1+z)s_1(z)}-\dfrac{1}{z}\\
&=\dfrac{((3-z)z-1-z)s_1(z)-z^2-z-1+3z}{z(1-3z+(1+z)s_1(z))}\\
&=-(1-z)^2\dfrac{s_1(z)+1}{z(1-3z+(1+z)s_1(z))}.
\end{split}
\]
Assume now that
\[
s(z)-\dfrac{1}{z} =O((1-z)^4).
\]
Then
\[
\dfrac{s_1(z)+1}{z(1-3z+(1+z)s_1(z))}=O((1-z)^2).
\]
In particular, $s_1$ satisfies the interpolation condition
\begin{equation}
\label{lumbago44}
s_1(1)=-1,\quad s_1^{\prime}(1)=0,
\end{equation}
and hence, by Lemma \ref{Richard-Lenoir}, $s_1(z)\equiv -1$. The corresponding
$s(z)=\frac{1}{z}$.

\section{The case where $s^\prime (1)=\alpha$}
\setcounter{equation}{0}
In this section we replace the condition $s^\prime(1)=1$ (see
\eqref{Crimee, ligne 7}) by $s^\prime(1)=\alpha$, where
$\alpha\in[0,1)$. Note that if $\alpha=0$, then it follows
immediately from Lemma \ref{Richard-Lenoir} that $s(z)\equiv 1$.
However, the fact that a Schur function $s$ satisfies the
interpolation conditions
\begin{equation}
s(1)=1,\quad s^\prime(1)=\alpha,\quad s^{(2)}(1)=s^{(3)}(1)=0
\label{lumbago}
\end{equation}
does not imply that $s(z)=\alpha z+1-\alpha$. Actually, for
each
$\alpha\in(0,1)$, there exists a sufficiently small $\beta>0$
such that the function
\[
s(z)=\alpha z+ 1-\alpha+\beta(1-z)^4
\]
is still a Schur function. For instance, for $\alpha=1/2$, the
following example is given in \cite{ELRS}:
\begin{equation}
s(z)=\frac{1+z}{2}+\frac{(z-1)^4}{20}.
\label{Daumesnil}
\end{equation}

To study the case $\alpha\in(0,1)$,
we first give, as in the previous section, a description of all Schur
functions $s$
such that
\begin{equation}
\label{parmentier}
s(z)=1+\alpha(z-1)+O((z-1)^4).
\end{equation}
We now have
\[
z_1=1,\quad k=1\quad{\rm and}\quad\tau_0=1,\quad \tau_1=\alpha.
\]
The matrix ${\mathbb P}$ given in \cite[(1.6), p. 3]{ADLW} reduces
to a number: ${\mathbb P}=\overline{\tau_0}\tau_1z_1=\alpha>0$.
The polynomial $p$ is now equal to
\[
p(z)=(1-z)\frac{1}{1-z}\frac{1}{\alpha}\frac{1}{2}=\frac{1}{2\alpha},\]
and we have
\[
\Theta_\alpha(z)=\begin{pmatrix}1&0\\0&1\end{pmatrix}-\frac{1}{2\alpha}
\frac{1+z}{1-z}\begin{pmatrix}1\\
1\end{pmatrix}\begin{pmatrix}1&-1\end{pmatrix}.\]
Thus
\[
\begin{split}
2\alpha(1-z)\Theta_\alpha(z)&=\begin{pmatrix}2\alpha(1-z)-(1+z)&1+z\\
-(1+z)&2\alpha(1-z)+1+z\end{pmatrix}\\
&=
\begin{pmatrix}
2\alpha-1-z(2\alpha+1)&1+z\\-1-z&2\alpha+1-z(1+2\alpha)\end{pmatrix}.
\end{split}
\]
A Schur function $s$ satisfies the interpolation conditions
\eqref{lumbago} if and only if it is of the form
\begin{equation}
\label{lumbago3} s(z)= \dfrac{(2\alpha-1-z(2\alpha+1))s_1(z)+1+z}
{2\alpha+1-z(2\alpha-1)-(1+z)s_1(z)},
\end{equation}
where, as in \eqref{lumbago1}, the parameter $s_1$ is any
Schur function satisfying \eqref{liminf}.\\

The parameter $s_1(z)=1-2\alpha$ satisfies \eqref{liminf}, and
corresponds to the solution
\[
s(z)=\dfrac{(2\alpha-1-z(2\alpha+1))(1-2\alpha)+1+z}
{2\alpha+1-z(2\alpha-1)-(1+z)(1-2\alpha)}=\alpha z+1-\alpha.
\]
Furthermore, with $s$ of the form \eqref{lumbago3},
\[
s(z)-(\alpha
z+1-\alpha)=\alpha(z-1)^2\dfrac{s_1(z)+2\alpha-1}{2\alpha+1-z(2\alpha-1)-(1+z)s_1(z)}.
\]
Assume that \eqref{parmentier} holds. Then
\[
\dfrac{s_1(z)+2\alpha-1}{2\alpha+1-z(2\alpha-1)-(1+z)s_1(z)}=O((1-z)^2).
\]
In particular, $s_1$ satisfies the interpolation conditions
\begin{equation}
\label{lumbago4} s_1(1)=1-2\alpha,\quad s_1^{\prime}(1)=0.
\end{equation}
However, as we have already mentioned above, these conditions do not
imply that $s_1(z)\equiv 1-2\alpha$.  For instance, for the
function defined by \eqref{Daumesnil}, we have
\[
s_1(z)=\frac{\dfrac{(z-1)^4}{10}}{2z+(1+z)\left(\dfrac{1+z}{2}
+\dfrac{(z-1)^4}{20}\right)},
\]
which indeed satisfies \eqref{lumbago4} with $\alpha = 1/2$, but is not
identically equal to $0$.\\

 At this stage we see
the difference between the cases $\alpha=1$ and $\alpha\in(0,1)$.
When $\alpha=1$, conditions \eqref{lumbago4} force $s_1(z)\equiv
-1$. When $\alpha\in(0,1)$, we need to impose more conditions on
$s$ in order to force $s_1(z)\equiv 1-2\alpha$. These conditions are
spelled out in the following theorem.

\begin{Tm}
Let $\alpha \in (0,1)$, assume that $s$ is a Schur function which admits
the representation
\[
s(z)=\alpha z+1-\alpha+O((1-z)^4)
\]
as $z$ tends nontangentially to $1$, and let $s_1$ be defined by
\eqref{lumbago3}. The following assertions are equivalent:\\
$(i)$ $s(z)=\alpha z +1-\alpha$.\\
$(ii)$ $s_1(z)\equiv 1-2\alpha$.\\
$(iii)$ For every $z\in{\mathbb D}$, it holds that $|s_1(z)|\le
|1-2\alpha|$.\\
$(iv)$ For every $z\in{\mathbb D}$, it holds that
\[
|(2\alpha+1-z(2\alpha-1))s(z)-1-z|\le|1-2\alpha|\cdot|
(2\alpha-1-z(2\alpha+1)-((2\alpha+1)s(z)-z(2\alpha-1))|.
\]
$(v)$ For every $z\in{\mathbb D}$, it holds that
\[
\frac{|1-s(z)|^2}{1-|s(z)|^2} <  \frac{\alpha}{1-\alpha}.
\]
\end{Tm}
{\bf Proof:} It is clear that
\[
(i)\iff(ii)\Rightarrow (iii)\iff (iv).\]
 Assume now that $(iv)$
is in force. If $\alpha =1/2$, then $(ii)$ and $(i)$ are clear.
If $\alpha \not=1/2$, consider the function $\sigma
=s_{1}/(1-2\alpha )$. By hypothesis, $\sigma $ is a Schur
function. If
its value at $z=1$ is $1$ and its angular derivative at $z=1$ is $0$, then $%
\sigma (z)\equiv 1$ by Lemma \ref{Richard-Lenoir}, that is, $(ii)$ holds.%
\newline

It remains to be shown that $(i)\iff (v)$. One direction is clear
by a direct computation. To prove the reverse assertion, we
consider the holomorphic function $f:{\mathbb{D}}\rightarrow \Pi
_{+}=\{w\in
\mathbb{C}
:\mathrm{Re}~w>0\}$ defined by
\begin{equation}
f(z)=C(s(z)),  \label{N5}
\end{equation}%
where $C(z)=\frac{1+z}{1-z}$ is the Cayley transform.
Calculations show that $f$ admits the representation
\begin{equation}
f(z)=\frac{1}{\alpha }\cdot \frac{1+z}{1-z}+\frac{1-\alpha }{\alpha }%
+r_{f}(z),  \label{N5'}
\end{equation}%
where
\begin{equation*}
r_{f}(z)=O((1-z)^{2}).
\end{equation*}%
In addition, condition $(v)$ implies that $\mathrm{Re}~f(z)>\frac{1-\alpha }{%
\alpha }$. So the function $f_{1}$, defined by
\begin{equation*}
f_{1}(z)=f(z)-\frac{1-\alpha }{\alpha },
\end{equation*}%
is of positive real part and
\begin{equation*}
(1-z)f_{1}\left( z\right) \rightarrow \frac{2}{\alpha },
\end{equation*}%
when $z\hat{\rightarrow}1$. Therefore, if we define a
self-mapping $h$\ of $\mathbb D $\ by $h=C^{-1}\left( f_{1}\right),
$\ then we have 
\begin{equation*}
(1-z)\frac{1+h(z)}{1-h(z)}\rightarrow \frac{2}{\alpha },
\end{equation*}%
when $z\hat{\rightarrow}1,$  which means that $h\left( 1\right) =1$\ and 
$%
h^{\prime }\left( 1\right) =\alpha .$\ Applying now the Julia-Wolff-Carath%
\'{e}odory theorem \cite{sarason94, RS} to the function $h$, we get that%
\begin{equation*}
\frac{|1-h(z)|^{2}}{1-|h(z)|^{2}}\leq \alpha
\frac{|1-z|^{2}}{1-|z|^{2}}
\end{equation*}%
or, equivalently,
\[
\mathrm{Re}~f_{1}(z)\geq \frac{1}{\alpha
}\mathrm{Re}~\frac{1+z}{1-z},
\]
because
\[
\mathrm{Re}~f_{1}(z) =\frac{1-|h(z)|^{2}}{|1-h(z)|^{2}}\quad{\rm
and}\quad \mathrm{Re}~\frac{1+z}{1-z}=\frac{1-|z|^{2}}{
|1-z|^{2}}.
\]
Thus we obtain that
\begin{equation*}
\mathrm{Re}~r_{f}(z)=\mathrm{Re}~\left( f_{1}(z)-\frac{1}{\alpha }\frac{1+z}{%
1-z}\right) \geq 0.
\end{equation*}%
Therefore the function $g$ defined by
\begin{equation*}
g\left( z\right) =\frac{1-r_{f}(z)}{1+r_{f}(z)}
\end{equation*}%
is a Schur function. In addition, $g\left(
1\right) =1$ and
the angular derivative $g^{\prime }\left( 1\right) =0.$ Therefore Lemma \ref%
{Richard-Lenoir} implies that
\begin{equation*}
g\left( z\right)  \equiv 1
\end{equation*}%
or
\begin{equation*}
r_{f}(z)\equiv 0.
\end{equation*}

Consequently, $f(z)=\frac{1}{\alpha }\cdot \frac{1+z}{1-z}+\frac{1-\alpha }{%
\alpha },$ or $s\left( z\right) =C^{-1}\left( f\left( z\right)
\right) =\alpha z+1-\alpha $, as asserted. \qed\mbox{}\\

\textbf{Remark:} Condition $(v)$ has a nice geometric
interpretation: the image $s({\mathbb{D}})$ of the open unit disk
$\mathbb{D}$ lies inside the horocycle
\begin{equation*}
{D(1, K)=\left\{z\in{\mathbb{D}} \, : \, \frac{|1-z|^2}{1-|z|^2} < K\right\}}%
,
\end{equation*}
of ``size" $K=\alpha/(1-\alpha)$, which is internally tangent to
the unit
circle ${\mathbb{T}}$ at the point $z=1$. See, for example, \cite[p. 121]{RS}%
\mbox{}\\

{\bf Acknowledgments:} the first author thanks the Earl Katz
family for endowing the chair which supports his research. The
second author was partially supported by the Fund for the
Promotion of Research at the Technion and by the Technion
President's Research Fund.
\bibliographystyle{plain}
\def\cprime{$'$} \def\cprime{$'$} \def\cprime{$'$} \def\cprime{$'$}
  \def\cprime{$'$}


\begin{thebibliography}{1}
\bibitem{Alpay01}
D.~Alpay,
\newblock {\em The {S}chur algorithm, reproducing kernel spaces and system
  theory},
\newblock SMF/AMS Texts and Monographs, Providence, RI, and Paris, 2001.

\bibitem{ADLW}
D.~Alpay, A.~Dijksma, H.~Langer, and G.~Wanjala,
\newblock {Basic boundary interpolation for generalized {S}chur functions and
  factorization of rational $J$--unitary matrix functions},
\newblock in D.~Alpay and I.~Gohberg, editors, {\em {Interpolation, Schur
  functions and moment problems}},
Operator Theory:
Advances and Applications,
volume 165, pages 1--29, Birkh{\" a}user
  Verlag, Basel, 2006.

\bibitem{adrs}
D.~Alpay, A.~Dijksma, J.~Rovnyak, and H.~de~Snoo,
\newblock {\em {Schur} functions, operator colligations, and reproducing kernel
  {P}ontryagin spaces},
Operator Theory:
Advances and Applications,
volume 96,
\newblock Birkh{\" a}user Verlag, Basel, 1997.
\bibitem{MR1242454}
D. M.~Burns and S. G. ~Krantz,
\newblock Rigidity of holomorphic mappings and a new {S}chwarz lemma at the
  boundary,
\newblock {\em J. Amer. Math. Soc.} 7 (1994), 661--676.

\bibitem{ELRS} M. Elin, M. Levenshtein, S. Reich, and D. Shoikhet,
Rigidity results for holomorphic mappings on the unit disk, in
{\em Complex and harmonic analysis}, pages 93-110, DEStech
Publications, Lancaster, PA, 2007.


\bibitem{kl1}
M. G. Kre{\u\i}n and H.~Langer,
\newblock {\"{U}}ber die verallgemeinerten {R}esolventen und die
  charakteristische {F}unktion eines isometrischen {O}perators im {R}aume ${\Pi
  _k}$,
\newblock in {\em Hilbert space operators and operator algebras (Proc.
Int.
  Conf., Tihany, 1970)}, pages 353--399,
\newblock Colloquia Math. Soc. J\'{a}nos Bolyai, 5,
North--Holland, Amsterdam, 1972.
\bibitem{RS} S. Reich and D. Shoikhet,
\newblock{\em Nonlinear semigroups, fixed points, and geometry
of domains in Banach spaces},
\newblock Imperial College Press,
London, 2005.
\bibitem{sarason94}
D.~Sarason,
\newblock {\em Sub--{H}ardy {H}ilbert spaces in the unit disk},
\newblock Wiley, {N}ew {Y}ork, 1994.

\end{thebibliography}

\end{document}